\documentclass[12pt]{article}
\usepackage{eurosym}
\usepackage{amsfonts}
\usepackage{amsmath,amssymb,amsthm}
\usepackage{geometry}
\usepackage[T1]{fontenc}
\usepackage[utf8]{inputenc}
\usepackage{lmodern}
\usepackage{amsmath,amssymb,mathtools}
\usepackage{xcolor}
\usepackage{subcaption}
\usepackage{tikz}

\usetikzlibrary{arrows.meta,calc,intersections}
\newtheorem{theorem}{Theorem}[section]

\newtheorem{remark}[theorem]{Remark}

\numberwithin{equation}{section}

\begin{document}

\title{Entire Large Solutions for Competitive Semilinear Elliptic Systems\\
with General Nonlinearities Satisfying Keller--Osserman Conditions}
\author{Dragos-Patru Covei\\[0.3cm]
{\small Department of Applied Mathematics}\\
{\small The Bucharest University of Economic Studies}\\
{\small Piata Romana, 1st district, Postal Code 010374, Bucharest, Romania}\\[0.2cm]
{\small E-mail: \texttt{dragos.covei@csie.ase.ro}}}
\date{}
\maketitle

\begin{abstract}
We study positive radial solutions of the weighted cooperative system
$\Delta u=p(|x|)f(u,v)$, $\Delta v=q(|x|)g(u,v)$ in $\mathbb{R}^{n}$,
$n\geq3$, when the Newtonian potentials of the nonnegative weights are
finite and the nonlinear source possesses a Keller--Osserman lower
envelope.  A threshold set of central data is introduced and its order and
compactness properties are used to construct an entire solution on its
interior boundary.  A decreasing nonlinear transform applied to
$W=u+v$ yields an explicit lower estimate in terms of the tail Newtonian
potential and proves that $W(r)\to\infty$.  Consequently at least one
component is large.  We also identify the additional hypothesis genuinely
needed for componentwise blow-up: diagonal Keller--Osserman integrability
alone is insufficient, whereas a two-sided balance of the weighted reaction
channels forces $u(r),v(r)\to\infty$.  Two independent proofs of this
transfer principle are given, one based on radial fluxes and one on
Newtonian potentials.  Finally, the PDE system is derived from an
isothermal two-species diffusion--reaction balance in a heterogeneous
porous catalyst, including nondimensional parameters and the interpretation
of the whole-space idealization.
\end{abstract}



\section{Introduction}

The study of \emph{positive entire large solutions} to nonlinear elliptic
systems is a fundamental area in the theory of partial differential equations, with profound connections to differential geometry, mathematical physics, reaction-diffusion systems, and population dynamics. Solutions that blow up at infinity arise naturally in various contexts: prescribed scalar curvature problems \cite{KO57,Os57}, conformal geometry, steady-state models of chemical reactions, and competitive or cooperative biological systems. 

The existence and qualitative behavior of such large solutions have been investigated extensively since the pioneering works of Keller \cite{KO57} and Osserman \cite{Os57} in the 1950s, who established the now-classical Keller--Osserman condition for scalar equations. This condition characterizes when a single semilinear equation $\Delta u = f(u)$ admits solutions that diverge to infinity at spatial infinity. Extensions to systems began with the work of \ \cite{Covei2017,Covei2017JMAA} on bounded solutions, and have seen significant developments in recent years, including the work of Li and Yang \cite{LiYang2015} on non-monotone systems, and more recently, Covei \cite{Covei2017,Covei2017JMAA} for competitive systems with power-type nonlinearities.

In this work, we focus on the radial
semilinear system%
\begin{equation}
\begin{cases}
\Delta u=p(r)f(u,v), \\ 
\Delta v=q(r)g(u,v),%
\end{cases}%
\quad r=|x|,\quad x\in \mathbb{R}^{n},\ n\geq 3,  \label{eq:system}
\end{equation}
in the entire space $\mathbb{R}^n$, where the weight functions $p,q$ are radially symmetric and satisfy appropriate integrability conditions, and the nonlinearities $f,g$ satisfy general Keller--Osserman-type growth conditions.

\medskip \noindent We assume:

\begin{itemize}
\item[(H1)] \textbf{Weights:} $p,q\geq 0$ are continuous, not both trivial, $%
\min (p,q)$ not compactly supported, with finite Green-type potentials%
\begin{equation*}
P(\infty )=\frac{1}{n-2}\int_{0}^{\infty }r\,p(r)\,dr<\infty ,\quad Q(\infty
)=\frac{1}{n-2}\int_{0}^{\infty }r\,q(r)\,dr<\infty ,
\end{equation*}%
where%
\begin{equation*}
P(r):=\int_{0}^{r}t^{1-n}\int_{0}^{t}s^{n-1}p(s)\,ds\,dt,\quad
Q(r):=\int_{0}^{r}t^{1-n}\int_{0}^{t}s^{n-1}q(s)\,ds\,dt,
\end{equation*}%
and 
\begin{equation*}
P(\infty )=\lim_{r\rightarrow \infty }P(r),Q(\infty )=\lim_{r\rightarrow
\infty }Q(r).
\end{equation*}

\item[(H2)] \textbf{Monotonicity:} $f,g$ are continuous on $[0,\infty )^{2}$%
, $C^{1}$ on $(0,\infty )^{2}$, positive on $(0,\infty )^{2}$, vanish at $%
(0,0)$, and are nondecreasing in each variable.

\item[(H3)] \textbf{Keller--Osserman envelope:} There exists a continuous,
nondecreasing $h\geq 0$, with $h(s)>0$ for $s>0$, such that for some small $%
\eta >0$,%
\begin{equation*}
f(s,t)+g(s,t)\geq h(s+t),\quad \forall s,t\geq \eta ,
\end{equation*}%
and the Keller--Osserman integral%
\begin{equation*}
H(\infty )=\lim_{s\rightarrow \infty }\int_{1}^{s}\frac{dt}{\sqrt{%
\int_{0}^{t}h(z)\,dz}}<\infty
\end{equation*}%
is finite.

\item[(H4)] \textbf{Balanced reaction channels (used only for componentwise
blow-up):} there is $\kappa\in(0,1]$ such that
\begin{equation}
\kappa\,q(r)g(s,t)\le p(r)f(s,t)\le
\kappa^{-1}q(r)g(s,t)
\label{eq:balance}
\end{equation}
for $r\ge0$ and $s,t\ge\eta$.  When componentwise blow-up is asserted we
also retain the diagonal integrability
\begin{equation}
\int_1^\infty\frac{dt}{f(t,t)}<\infty,\qquad
\int_1^\infty\frac{dt}{g(t,t)}<\infty .
\label{kelleroser}
\end{equation}
\end{itemize}

\medskip\noindent\textbf{Main contributions.}
Our principal aim is to extend the existence theory of Li-Yang~%
\cite{LiYang2015}, originally established for power-type nonlinearities, to a broad class
of nonlinearities satisfying Keller--Osserman-type growth conditions. This extension is both natural and necessary:

\begin{itemize}
\item[(i)] \emph{Generality:} Many physical and geometric problems lead to nonlinearities that do not have pure power-law form. Our framework accommodates these more general growth conditions while preserving the essential analytical structure.

\item[(ii)] \emph{Unification:} By working with a general envelope function $h$ satisfying the Keller--Osserman integral condition, we provide a unified treatment of critical and supercritical growth regimes, recovering as special cases not only Covei's~\cite{Covei2009} results but also classical results of Keller~\cite{KO57} and Osserman~\cite{Os57}.

\item[(iii)] \emph{Completeness:} Our results complement the bounded-solution theories developed by Covei~\cite{Covei2017,Covei2017JMAA}, Li and Yang~\cite{LiYang2015}, and Santos et al.~\cite{Santos2017}, providing a complete picture of the existence landscape for system \eqref{eq:system}.
\end{itemize}

\medskip\noindent\textbf{Analytical framework.}
We call $(u,v)$ a \emph{positive entire solution} if it
is $C^2$, radial, and strictly positive on $[0,\infty)$. A solution is called \emph{large} if at least one component diverges to infinity as $r \to \infty$. The set of \emph{%
central values}
\begin{equation*}
\mathcal{G}=\{(\alpha ,\beta )\in \mathbb{R}^{+}\times \mathbb{R}^{+}\mid
\exists \ \text{entire solution of \eqref{eq:system} with}\
(u(0),v(0))=(\alpha ,\beta )\}
\end{equation*}%
encodes the initial data leading to such solutions. A key insight of our work is that $\mathcal{G}
$ has a rich geometric structure: it is downward-closed (in the sense that smaller initial values also lead to entire solutions), certain regularized subsets $\mathcal{G}_\delta$ are compact, and solutions corresponding to boundary points of $\mathcal{G}_\delta$ are necessarily large---at least one component diverges as $%
r\rightarrow \infty $.

\medskip\noindent\textbf{Methodology.}
Our approach adapts and extends Covei's barrier and radial integration
techniques~\cite{Covei2009,Covei2017,Covei2017JMAA} to the setting of general nonlinearities. The key steps are:
\begin{itemize}
\item \emph{Monotone iteration:} We construct solutions via a monotone iterative scheme, with upper bounds provided by scalar supersolutions.
\item \emph{Scalar reduction:} By exploiting the monotonicity of $f$ and $g$, we reduce the coupled system to a scalar differential inequality for $W = u + v$.
\item \emph{Keller--Osserman transform:} We apply a nonlinear transformation $H$ to $W$, converting the nonlinear problem into a linear inequality amenable to explicit integration.
\item \emph{Blow-up dichotomy:} Using barrier arguments and the Hopf lemma, we establish that solutions with boundary initial data must blow up.
\end{itemize}
This approach not only recovers known power-type results but also provides explicit quantitative bounds on the blow-up rate.

\medskip \noindent \textbf{Organization.} Section~\ref{2} states
the main results. Section~\ref{3} contains the complete proofs, built on monotone
iteration, Keller--Osserman estimates, and barrier arguments. Section~\ref{aplicatie} outlines the physical motivation that leads to the mathematical model studied in this work.  Section~\ref{5} summarizes the conclusions.

\section{Main Results \label{2}}

We begin with the following property of the set $\mathcal{G}$.

\begin{theorem}
\label{lem:G-nonempty} Let $\alpha_{0} > 0$, $\beta_{0} > 0$. Then $\mathcal{%
G}$ is nonempty and has the property that if $(\alpha_{0}, \beta_{0}) \in 
\mathcal{G}$ and $0 < \alpha < \alpha_{0}$, $0 < \beta < \beta_{0}$, then $%
(\alpha, \beta) \in \mathcal{G}$.
\end{theorem}

For $\delta >0$ sufficiently small, define 
\begin{equation*}
\mathcal{G}_{\delta }:=\left\{ (\alpha ,\beta )\in \mathcal{G}\ \middle|\
\min \{\alpha ,\beta \}\geq \delta \right\} .
\end{equation*}%
We have the following result.

\begin{theorem}
\label{lem:Gphi} The set $\mathcal{G}_{\delta}$ is both closed and bounded
in $(0,\infty) \times (0,\infty)$.
\end{theorem}

Next, we define the \emph{edge set} 
\begin{equation*}
\mathcal{E}(\mathcal{G}_{\delta }):=\left\{ (\alpha ,\beta )\in \partial 
\mathcal{G}_{\delta }\ \middle|\ \min \{\alpha ,\beta \}>\delta \right\} .
\end{equation*}%
The geometric relationship between these sets is depicted in Figure 1, where 
$\mathcal{G}$ is shown in green, $\mathcal{G}_{\delta }$ in yellow, and the
interior boundary arc $\mathcal{E}(\mathcal{G}_{\delta })$ in blue.

\begin{figure}[h!]
\centering
\begin{tikzpicture}[scale=1.05, >=Stealth, line cap=round, line join=round]
  % Axes
  \draw[->] (-0.1,0) -- (8.1,0) node[below right] {$\alpha$};
  \draw[->] (0,-0.1) -- (0,6.1) node[above left] {$\beta$};

  % Parameter delta
  \def\del{1.00}

  % Fill G (green)
  \fill[green!30] (0,0) -- (0, 5.5) .. controls (2.0, 5.0) and (5.0, 2.0) .. (5.5, 0) -- cycle;
  
  % Outline of G boundary
  \draw[green!50!black, very thick] (0, 5.5) .. controls (2.0, 5.0) and (5.0, 2.0) .. (5.5, 0);

  % Guide lines for alpha = delta and beta = delta
  \draw[gray!70, dashed] (\del,0) -- (\del,6.1) node[above, black] {\small $\alpha=\delta$};
  \draw[gray!70, dashed] (0,\del) -- (8.1,\del) node[right, black] {\small $\beta=\delta$};

  % Fill G_delta (yellow)
  \begin{scope}
    \clip (\del,\del) rectangle (8.1,6.1);
    \fill[yellow!50] (0,0) -- (0, 5.5) .. controls (2.0, 5.0) and (5.0, 2.0) .. (5.5, 0) -- cycle;
  \end{scope}

  % E(G_delta) = boundary arc with alpha>delta, beta>delta (blue)
  \begin{scope}
    \clip (\del,\del) rectangle (8.1,6.1);
    \draw[blue!80!black, ultra thick] (0, 5.5) .. controls (2.0, 5.0) and (5.0, 2.0) .. (5.5, 0);
  \end{scope}

  % Labels
  \node[green!50!black, fill=white, inner sep=1pt, rounded corners=2pt] at (0.5, 0.5) {$\mathcal{G}$};
  \node[black, fill=white, inner sep=1pt, rounded corners=2pt] at (2.2, 2.2) {$\mathcal{G}_{\delta}$};
  \node[blue!80!black, fill=white, inner sep=1pt, rounded corners=2pt] at (3.6, 3.6) {$\mathcal{E}(\mathcal{G}_{\delta})$};

  % Tick marks
  \foreach \x in {1,2,3,4,5,6,7} \draw (\x,0) -- ++(0,0.06);
  \foreach \y in {1,2,3,4,5,6} \draw (0,\y) -- ++(0.06,0);
\end{tikzpicture}
\caption{Existence region in the $(\protect\alpha,\protect\beta)$-plane: $\mathcal{G}$
(green), $\mathcal{G}_{\protect\delta}$ (yellow), and $\mathcal{E}(\mathcal{G}_{\protect\delta})$ (blue).}
\end{figure}
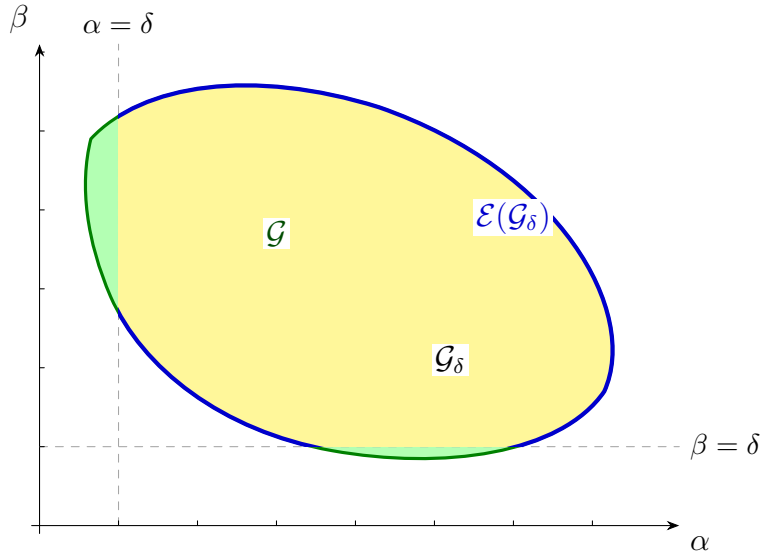

We can now state our main existence theorem for entire large solutions of %
\eqref{eq:system}.

\begin{theorem}
\label{thm:main} There exists a positive radial solution $(u,v)$ of %
\eqref{eq:system} with central value $(u(0),v(0))\in \mathcal{E}(\mathcal{G}%
_{\delta })$ such that at least one of $u$ or $v$ diverges, i.e., 
\begin{equation}
\lim_{r\rightarrow \infty }u(r)=\infty \quad \text{or}\quad
\lim_{r\rightarrow \infty }v(r)=\infty .  \label{el}
\end{equation}%
Moreover, if $H(s):=\int_{s}^{\infty }\frac{dt}{f(t,t)+g\left( t,t\right) }$
and $H^{-1}$ denotes its inverse, then for all sufficiently large $r$, 
\begin{equation}
u(r)+v(r)\ \geq \ H^{-1}\!\left( \int_{r}^{\infty }t^{1-n}\int_{0}^{t}s^{n-1}%
\big(p(s)+q(s)\big)\,ds\,dt\right) .  \label{eq:lower-bound}
\end{equation}%
If, in addition, the balanced-channel hypothesis\/ {\rm(H4)} holds, then
\[
\lim_{r\to\infty}u(r)=\lim_{r\to\infty}v(r)=\infty .
\]
\end{theorem}

\begin{remark}[Relation to Covei's result]
Both Theorem~\ref{thm:main} and \cite{Covei2017,Covei2017JMAA} study radial
elliptic systems of the form~\eqref{eq:system}. The two results are
therefore \emph{compatible}: they describe different behaviours (bounded
vs.\ large) in disjoint regions of the central value space, even for the
same coefficients.

In \cite{Covei2017,Covei2017JMAA}, the focus is on positive \emph{bounded
entire} radial solutions, under a strict hierarchy between the finite
numbers $P(\infty )$, $Q(\infty )$ and $H(\infty )$.

Theorem~\ref{thm:main} asserts that for $(u(0),v(0))\in \mathcal{E}(\mathcal{%
G}_{\delta })$, any entire solution is large, while the results of \cite%
{Covei2017,Covei2017JMAA} produces bounded entire solutions for $(u(0),v(0))$ in the
interior of $\mathcal{G}$.
\end{remark}

\begin{remark}
The sets $\mathcal{G}$, $\mathcal{G}_{\delta }$, and $\mathcal{E}(\mathcal{G}%
_{\delta })$ can be thought of as nested regions: $\mathcal{G}$ contains all
starting values that lead to positive solutions, $\mathcal{G}_{\delta }$ is
a smaller core within it, and $\mathcal{E}(\mathcal{G}_{\delta })$ marks the
area where the solutions grow without bound.
\end{remark}

\section{Proof of the Theorems \label{3}}

We recall that any radial $C^{2}$ solution $(u,v)$ of \eqref{eq:system}
satisfies the equivalent Volterra-type integral system 
\begin{align}
u(r) &= \alpha_{0} + \int_{0}^{r} t^{1-n} \int_{0}^{t} s^{n-1} p(s) \, f\big(%
u(s), v(s)\big) \, ds \, dt,  \label{eq:int-u} \\
v(r) &= \beta_{0} + \int_{0}^{r} t^{1-n} \int_{0}^{t} s^{n-1} q(s) \, g\big(%
u(s), v(s)\big) \, ds \, dt,  \label{eq:int-v}
\end{align}
for all $r \geq 0$, where $\alpha_{0} = u(0)$ and $\beta_{0} = v(0)$ are the
central values. For $\alpha_{0}, \beta_{0} \in (0,\infty)$, any solution of %
\eqref{eq:int-u}--\eqref{eq:int-v} is also a solution of \eqref{eq:system},
though it may not satisfy the large-solution condition \eqref{el}.

\subsection{Proof of Theorem~\protect\ref{lem:G-nonempty}}

\begin{proof}
The proof is organized into six logical steps. We adapt the monotone iteration method of \cite{Covei2017,Covei2017JMAA,LiYang2015} to the
present setting, where $f$ and $g$ satisfy the general hypotheses (H2)--(H3).

\medskip\noindent \textbf{Step 1: Iterative scheme.}
Given $(\alpha _{0},\beta _{0})\in (0,\infty )\times (0,\infty )$, define%
\begin{equation*}
u_{0}(r)\equiv \alpha _{0},\quad v_{0}(r)\equiv \beta _{0},
\end{equation*}%
and for $k\geq 1$ set 
\begin{align*}
u_{k}(r)& =\alpha _{0}+\int_{0}^{r}t^{1-n}\int_{0}^{t}s^{n-1}p(s)\,f\big(%
u_{k-1}(s),v_{k-1}(s)\big)\,ds\,dt, \\
v_{k}(r)& =\beta _{0}+\int_{0}^{r}t^{1-n}\int_{0}^{t}s^{n-1}q(s)\,g\big(%
u_{k-1}(s),v_{k-1}(s)\big)\,ds\,dt.
\end{align*}

\medskip\noindent \textbf{Step 2: Monotonicity.}
By the positivity of $p,q$ and the monotonicity of $f,g$ in each variable
(assumption (H2)), we have%
\begin{equation*}
u_{k}(r)\geq u_{k-1}(r),\quad v_{k}(r)\geq v_{k-1}(r),\quad \forall r\geq
0,\ \forall k\geq 1.
\end{equation*}%
Thus, for each fixed $r$, both sequences $\{u_{k}(r)\}$ and $\{v_{k}(r)\}$
are nondecreasing.

\medskip\noindent \textbf{Step 3: Uniform upper bound via a scalar supersolution.}
\emph{Construction of $w$.} The Keller--Osserman condition (H3) ensures the
existence of a scalar function $w\in C^{2}([0,\infty ))$ satisfying 
\begin{equation}
\Delta w=(p+q)\big(f(w,w)+g(w,w)\big),\quad w(0)>\alpha _{0}+\beta
_{0},\quad \lim_{r\rightarrow \infty }w(r)=\infty .  \label{super}
\end{equation}%
To construct $w$ rigorously, consider the scalar radial ODE%
\begin{equation}
W^{\prime \prime }(r)+\frac{n-1}{r}W^{\prime }(r)=(p(r)+q(r))\,S(W(r)),
\label{eq:W-ODE}
\end{equation}%
where $S(s):=f(s,s)+g(s,s)$. By (H2)--(H3), $S$ is continuous,
nondecreasing, positive for $s>0$, vanishes at $s=0$, and satisfies the Keller--Osserman
integral condition%
\begin{equation}
\int_{1}^{\infty }\frac{ds}{\sqrt{\int_{0}^{s}S(\tau )\,d\tau }}<\infty .
\label{eq:KO-S}
\end{equation}%
By classical existence theory for radial ODEs (see \cite{KO57,Os57}), for any initial value $W(0) = W_0 > 0$,
there exists a maximal solution $W$ of \eqref{eq:W-ODE} on an interval $[0, R_{W_0})$ with $R_{W_0} \in (0,\infty]$. 
The Keller--Osserman condition \eqref{eq:KO-S} guarantees that if $W_0$ is chosen sufficiently large, then $R_{W_0} = \infty$ and $W(r) \to \infty$ as $r \to \infty$. 
Specifically, choosing $W_0 > \alpha_0 + \beta_0$, we obtain a global solution $W$ satisfying
$W(0) = W_0 > \alpha_0 + \beta_0$ and $\lim_{r \to \infty} W(r) = \infty$. Setting $w := W$ yields \eqref{super}.

\emph{Significance of $w$.} The function $w$ is a \emph{scalar supersolution}
for the system: it bounds the sum $u+v$ from above and grows unboundedly at
infinity. This is crucial because:

\begin{itemize}
\item It provides a uniform \emph{a priori} bound on the iterates $%
u_{k}+v_{k}$, independent of $k$, on every compact interval.

\item It allows the use of comparison principles to control the nonlinear
coupling in the system by a single scalar equation.

\item It ensures that the iteration cannot blow
up prematurely on finite intervals.
\end{itemize}

By the comparison principle for radial integral equations (which follows from the maximum principle for the Laplacian), 
since $u_0(r) + v_0(r) = \alpha_0 + \beta_0 < w(0)$ and
\begin{equation*}
\Delta(u_{k}+v_{k}) = p(r)f(u_k,v_k) + q(r)g(u_k,v_k) \leq (p(r)+q(r)) \big(f(u_k+v_k,u_k+v_k) + g(u_k+v_k,u_k+v_k)\big),
\end{equation*}
by monotonicity of $f,g$, we deduce by induction that
\begin{equation*}
u_{k}(r)+v_{k}(r)\leq w(r),\quad \forall r\geq 0,\ \forall k\geq 0.
\end{equation*}

\medskip\noindent \textbf{Step 4: Convergence of the iteration.}
For each fixed $r\geq 0$, the sequences $\{u_{k}(r)\}$ and $\{v_{k}(r)\}$
are monotone increasing (by Step 2) and bounded above by $w(r)$ (by Step 3). Thus, the pointwise
limits%
\begin{equation*}
u(r):=\lim_{k\rightarrow \infty }u_{k}(r),\quad v(r):=\lim_{k\rightarrow
\infty }v_{k}(r)
\end{equation*}%
exist for all $r\geq 0$ by monotone convergence of real sequences.

\emph{Justification via the Monotone Convergence Theorem (MCT).} 
To verify that $(u,v)$ satisfies the integral system \eqref{eq:int-u}--\eqref{eq:int-v}, we apply Lebesgue's MCT.
The MCT states: If $\{\phi_{k}\}$ is
a sequence of nonnegative measurable functions such that $\phi_{k}(x)\uparrow
\phi(x)$ pointwise almost everywhere, then%
\begin{equation*}
\lim_{k\rightarrow \infty }\int \phi_{k}\,d\mu = \int \lim_{k\rightarrow \infty }\phi_{k}\,d\mu = \int \phi\,d\mu.
\end{equation*}%
In our setting, fix $r > 0$ and consider the innermost integral in \eqref{eq:int-u}:
\begin{equation*}
I_k(t) := \int_0^t s^{n-1} p(s) f(u_{k-1}(s), v_{k-1}(s))\,ds.
\end{equation*}
Since $u_{k-1}(s) \uparrow u(s)$ and $v_{k-1}(s) \uparrow v(s)$ pointwise, and $f$ is nondecreasing and continuous, we have
\begin{equation*}
p(s) f(u_{k-1}(s), v_{k-1}(s)) \uparrow p(s) f(u(s), v(s)) \quad \text{pointwise in } s.
\end{equation*}
The integrands are nonnegative. By MCT,
\begin{equation*}
\lim_{k\to\infty} I_k(t) = \int_0^t s^{n-1} p(s) f(u(s), v(s))\,ds.
\end{equation*}
Moreover, $I_k(t) \uparrow$ and is bounded by
\begin{equation*}
I_k(t) \leq \int_0^t s^{n-1} (p(s)+q(s)) S(w(s))\,ds =: M(t) < \infty
\end{equation*}
for each $t < \infty$, by construction of $w$.
Applying MCT again to the outer integral in \eqref{eq:int-u}, we obtain
\begin{equation*}
\lim_{k\to\infty} u_k(r) = \alpha_0 + \int_0^r t^{1-n} \left(\lim_{k\to\infty} I_k(t)\right) dt = \alpha_0 + \int_0^r t^{1-n} \int_0^t s^{n-1} p(s) f(u(s),v(s))\,ds\,dt.
\end{equation*}
Thus $u = \lim_{k\to\infty} u_k$ satisfies \eqref{eq:int-u}. The same argument applies to $v$. By standard elliptic regularity, $(u,v) \in C^2([0,\infty)) \times C^2([0,\infty))$.

\medskip\noindent \textbf{Step 5: Conclusion.}
The limit pair $(u,v)$ is a positive $C^{2}$ radial solution of %
\eqref{eq:system} with $(u(0),v(0))=(\alpha _{0},\beta _{0})$. Thus $(\alpha
_{0},\beta _{0})\in \mathcal{G}$, proving that $\mathcal{G}$ is nonempty.

\medskip\noindent \textbf{Step 6: Downward closure property.}
If $(\alpha _{0},\beta _{0})\in \mathcal{G}$ and $0<\alpha <\alpha _{0}$, $%
0<\beta <\beta _{0}$, then starting the iteration from $(\alpha ,\beta )$
yields smaller initial data. By the same monotonicity and boundedness
arguments (Steps 1--4), the iteration converges to a positive entire solution with initial values $(\alpha, \beta)$, hence $%
(\alpha ,\beta )\in \mathcal{G}$. This completes the proof of Theorem~\ref{lem:G-nonempty}.
\end{proof}

\subsection{Proof of Theorem~\protect\ref{lem:Gphi}}

\begin{proof}
We divide the proof into two logically independent parts: (i) \emph{%
boundedness} of $\mathcal{G}_{\delta }$, and (ii) \emph{closedness} of $%
\mathcal{G}_{\delta }$.

\medskip\noindent\textbf{Part I: Boundedness of $\mathcal{G}_{\delta}$.}

\emph{Step 0: Contradiction setup.} Suppose, for contradiction, that $%
\mathcal{G}_{\delta }$ is unbounded in $(0,\infty )\times (0,\infty )$. From
Theorem~\ref{lem:G-nonempty}, we know that $\mathcal{G}$ is \emph{downward
closed} in each coordinate: if $(\alpha _{0},\beta _{0})\in \mathcal{G}$ and 
$0<\alpha <\alpha _{0}$, $0<\beta <\beta _{0}$, then $(\alpha ,\beta )\in 
\mathcal{G}$.

In particular, if $(\alpha _{0},\beta _{0})\in \mathcal{G}_{\delta }$, then
for any $(\alpha ,\beta )$ with%
\begin{equation*}
\delta \leq \alpha \leq \alpha _{0},\quad \delta \leq \beta \leq \beta _{0},
\end{equation*}%
we also have $(\alpha ,\beta )\in \mathcal{G}_{\delta }$. This implies that
if $\mathcal{G}_{\delta }$ were unbounded, it would necessarily contain an
unbounded subset of the form%
\begin{equation*}
\lbrack \delta ,\infty )\times \{\delta \}\quad \text{or}\quad \{\delta
\}\times \lbrack \delta ,\infty ).
\end{equation*}%
Without loss of generality, assume%
\begin{equation*}
\{\delta \}\times \lbrack \delta ,\infty )\subseteq \mathcal{G}_{\delta }.
\end{equation*}%
\medskip \noindent \emph{Step 1: Choice of a scalar comparison function.}
Define%
\begin{equation*}
\widetilde{p}(r):=\min \{p(r),\,q(r)\}.
\end{equation*}%
By assumption (H1), $\widetilde{p}$ is continuous, radial, and not compactly
supported. Hence there exists $R>0$ such that $\widetilde{p}(R)>0$. From 
\cite{Covei2009} (and the Keller--Osserman framework), there exists a \emph{%
large} nonnegative radial solution $w$ of%
\begin{equation*}
\Delta w=\widetilde{p}(r)\,h(w)\quad \text{in }[0,R),\quad w(0)=\delta
,\quad \lim_{r\rightarrow R^{-}}w(r)=\infty .
\end{equation*}%
Here \emph{large} means that $w$ blows up at the finite radius $R$.

Significance of $w$:

The function $w$ is a \emph{scalar supersolution} for the coupled system: it
satisfies the same type of equation as $u+v$ but with the minimal weight $%
\widetilde{p}$ and the scalar nonlinearity $h$. It will serve as a \emph{%
barrier function} in the comparison argument: if $u+v$ starts above $w$ at $%
r=0$, the maximum principle forces $u+v$ to stay above $w$ until blow-up,
which is impossible if $u+v$ is smooth on $[0,R)$.

\medskip \noindent \emph{Step 2: Construction of a candidate solution in $%
\mathcal{G}_{\delta }$.} Since $\{\delta \}\times \lbrack \delta ,\infty
)\subseteq \mathcal{G}_{\delta }$, pick $\alpha _{0}>\delta $ and consider
the initial data%
\begin{equation*}
u(0)=\alpha _{0}>w(0)=\delta ,\quad v(0)=\delta .
\end{equation*}%
By Theorem~\ref{lem:G-nonempty}, there exists a positive entire solution $%
(u,v)$ of \eqref{eq:int-u}--\eqref{eq:int-v} with these initial values: 
\begin{align*}
u(r)& =\alpha
_{0}+\int_{0}^{r}t^{1-n}\int_{0}^{t}s^{n-1}p(s)\,f(u(s),v(s))\,ds\,dt, \\
v(r)& =\delta
+\int_{0}^{r}t^{1-n}\int_{0}^{t}s^{n-1}q(s)\,g(u(s),v(s))\,ds\,dt,
\end{align*}%
for $0\leq r<\infty $.

\medskip \noindent \emph{Step 3: Differential inequality for $u+v$.} Since $%
u(r)\geq \alpha _{0}>\delta $ and $v(r)\geq \delta $ for all $r\geq 0$, we
have 
\begin{equation*}
u(r)+v(r)\geq \alpha _{0}+\delta >2\delta .
\end{equation*}%
By assumption (H3), $f(s,t)+g(s,t)\geq h(s+t)$ for $s,t\geq \eta $ (and here 
$\eta \leq \delta $), so%
\begin{equation*}
\Delta (u+v)=p(r)f(u,v)+q(r)g(u,v)\geq \widetilde{p}(r)\,h(u+v)\quad \text{%
in }[0,R).
\end{equation*}%
\medskip \noindent \emph{Step 4: Comparison with $w$ and contradiction.} The
function $w$ satisfies $\Delta w=\widetilde{p}(r)\,h(w)$ in $[0,R)$, $%
w(0)=\delta $, and blows up at $r=R$. By the strong maximum principle
applied to $u+v-w$, if $u(0)+v(0)>w(0)$, then $u+v>w$ near $r=0$. However,
since $w$ blows up at $R$ while $u+v$ is $C^{2}$ on $[0,R)$, the comparison
principle forces $u+v\leq w$ in $[0,R)$ a contradiction with $u(0)+v(0)>w(0)$%
.

Therefore, $\mathcal{G}_{\delta }$ cannot be unbounded; it is bounded.

\medskip\noindent\textbf{Part II: Closedness of $\mathcal{G}_{\delta}$.}

\emph{Step 1: Sequential closure argument.} Let $\{(\alpha
_{m},\beta _{m})\}\subset \mathcal{G}_{\delta }$ be a sequence converging to 
$(\alpha ,\beta )\in \overline{\mathcal{G}_{\delta }}$. For each $m$,
Theorem~\ref{lem:G-nonempty} guarantees the existence of a positive entire
solution $(u_{m},v_{m})$ of \eqref{eq:int-u}--\eqref{eq:int-v} with%
\begin{equation*}
u_{m}(0)=\alpha _{m},\quad v_{m}(0)=\beta _{m}.
\end{equation*}%

\emph{Step 2: Uniform bounds and equicontinuity.} From
the monotone iteration scheme in Theorem~\ref{lem:G-nonempty} and the
Keller--Osserman supersolution $w$ from \eqref{super} with $w(0)>\sup_m(\alpha
_{m}+\beta _{m})$, we deduce that
\begin{equation*}
u_m(r) + v_m(r) \leq w(r) \quad \forall r \geq 0, \, \forall m \geq 1.
\end{equation*}
Hence the family $\{u_{m}\}$ and $\{v_m\}$ is uniformly
bounded on each compact interval $[0,R]$. 

To establish equicontinuity, we use the integral representation \eqref{eq:int-u}--\eqref{eq:int-v}.
For $r_1, r_2 \in [0,R]$ with $r_1 < r_2$, we have
\begin{align*}
|u_m(r_2) - u_m(r_1)| &\leq \int_{r_1}^{r_2} t^{1-n} \int_0^t s^{n-1} p(s) f(u_m(s), v_m(s))\,ds\,dt \\
&\leq \int_{r_1}^{r_2} t^{1-n} \int_0^t s^{n-1} (p(s)+q(s)) S(w(s))\,ds\,dt,
\end{align*}
where $S(s) = f(s,s) + g(s,s)$. Since the right-hand side is independent of $m$ and tends to zero as $r_2 - r_1 \to 0$ (by absolute continuity of the integral), the family $\{u_m\}$ is equicontinuous on $[0,R]$. The same argument applies to $\{v_m\}$.

\emph{Step 3: Compactness and passage to the limit.} By
the Arzel\`{a}--Ascoli theorem, the sequence $\{(u_m, v_m)\}$ is relatively compact in $C([0,R]) \times C([0,R])$ for each $R > 0$.
By a standard diagonal argument, there exists a subsequence $(u_{m_j}, v_{m_j})$ converging
uniformly on every compact interval $[0,R]$ to a limit $(u,v) \in C([0,\infty)) \times C([0,\infty))$. 

To verify that $(u,v)$ satisfies \eqref{eq:int-u}--\eqref{eq:int-v}, fix $r > 0$ and observe that
the integrands in \eqref{eq:int-u} converge:
\begin{equation*}
s^{n-1} p(s) f(u_{m_j}(s), v_{m_j}(s)) \to s^{n-1} p(s) f(u(s), v(s))
\end{equation*}
pointwise in $s \in [0,r]$ by continuity of $f$, and are uniformly bounded by $s^{n-1}(p(s)+q(s))S(w(s)) \in L^1([0,r])$.
By the dominated convergence theorem,
\begin{equation*}
\lim_{j \to \infty} \int_0^t s^{n-1} p(s) f(u_{m_j}(s), v_{m_j}(s))\,ds = \int_0^t s^{n-1} p(s) f(u(s), v(s))\,ds
\end{equation*}
uniformly in $t \in [0,r]$. Applying dominated convergence again to the outer integral in \eqref{eq:int-u}, we conclude that $(u,v)$ satisfies \eqref{eq:int-u}--\eqref{eq:int-v}, hence is a positive entire solution with%
\begin{equation*}
u(0)=\lim_{j\to\infty} u_{m_j}(0) = \lim_{j\to\infty}\alpha_{m_j} = \alpha ,\quad v(0)=\beta .
\end{equation*}%
Since $\min \{\alpha ,\beta \}\geq \delta $ by construction (as the limit of a sequence in $\mathcal{G}_\delta$), we have $%
(\alpha ,\beta )\in \mathcal{G}_{\delta }$.

\emph{Conclusion:} $\mathcal{G}_{\delta }$ is both bounded and closed, completing the proof of Theorem~\ref{lem:Gphi}.
\end{proof}

\subsection{Proof of Theorem~\protect\ref{thm:main}}

\begin{proof}
We structure the proof into seven distinct and logically self-contained steps.

\medskip\noindent\textbf{Standing assumptions and notation.}
We work under hypotheses (H1)--(H3) and consider the radial system %
\eqref{eq:system} in its equivalent integral (Green) form \eqref{eq:int-u}--%
\eqref{eq:int-v}, with prescribed central values%
\begin{equation*}
(u(0),v(0))=(\alpha ,\beta )\in \mathcal{E}(\mathcal{G}_{\delta }).
\end{equation*}%
We introduce the \emph{Keller--Osserman transform}%
\begin{equation*}
H(s):=\int_{s}^{\infty }\frac{dt}{f(t,t)+g(t,t)},\quad s>0.
\end{equation*}%
By construction:

\begin{itemize}
\item $H$ is strictly decreasing, with%
\begin{equation*}
H^{\prime }(s)=-\frac{1}{f(s,s)+g(s,s)}<0;
\end{equation*}

\item $H^{\prime \prime }(s)\geq 0$ since $s\mapsto f(s,s)+g(s,s)$ is
nondecreasing by (H2). Explicitly,
\begin{equation*}
H^{\prime\prime}(s) = \frac{d}{ds}\left(-\frac{1}{f(s,s)+g(s,s)}\right) = \frac{\frac{\partial}{\partial s}[f(s,s)+g(s,s)]}{[f(s,s)+g(s,s)]^2} \geq 0.
\end{equation*}
\end{itemize}

\medskip\noindent\textbf{Step 1: Perturbation of initial data and blow-up for approximants.}
From Theorem~\ref{lem:G-nonempty}, $\mathcal{G}$ is nonempty and downward
closed in each coordinate. From Theorem~\ref{lem:Gphi}, $\mathcal{G}_{\delta
}$ is closed and bounded for $\delta >0$ sufficiently small.

Let $(\alpha ,\beta )\in \mathcal{E}(\mathcal{G}_{\delta })$. For each $%
k\geq 1$, consider the perturbed initial data%
\begin{equation*}
\left( \alpha +\frac{1}{k},\ \beta +\frac{1}{k}\right) .
\end{equation*}%
If this pair belonged to $\mathcal{G}$, then by downward closure and the
definition of $\mathcal{G}_{\delta }$, $(\alpha ,\beta )$ would lie in the 
\emph{interior} of $\mathcal{G}_{\delta }$, contradicting $(\alpha ,\beta
)\in \mathcal{E}(\mathcal{G}_{\delta })$. Hence%
\begin{equation*}
\left( \alpha +\frac{1}{k},\ \beta +\frac{1}{k}\right) \notin \mathcal{G}.
\end{equation*}%
Therefore, the solution $(u_{k},v_{k})$ of 
\begin{align*}
u_{k}(r)& =\alpha +\frac{1}{k}+\int_{0}^{r}t^{1-n}%
\int_{0}^{t}s^{n-1}p(s)f(u_{k}(s),v_{k}(s))\,ds\,dt, \\
v_{k}(r)& =\beta +\frac{1}{k}+\int_{0}^{r}t^{1-n}%
\int_{0}^{t}s^{n-1}q(s)g(u_{k}(s),v_{k}(s))\,ds\,dt,
\end{align*}%
must blow up in finite radius: there exists $R_{k}>0$ such that%
\begin{equation*}
\lim_{r\rightarrow R_{k}^{-}}u_{k}(r)=\infty \quad \text{or}\quad
\lim_{r\rightarrow R_{k}^{-}}v_{k}(r)=\infty .
\end{equation*}%
In either case,%
\begin{equation}
\lim_{r\rightarrow R_{k}^{-}}\big(u_{k}(r)+v_{k}(r)\big)=\infty .  \label{le}
\end{equation}%
Moreover, by monotonicity in $k$ of the blow-up radii,%
\begin{equation*}
R_{1}\leq R_{2}\leq \dots \leq R_{k}\leq \dots ,
\end{equation*}%
and we set%
\begin{equation*}
R:=\lim_{k\rightarrow \infty }R_{k}\in \lbrack R_{1},\infty ].
\end{equation*}

\medskip\noindent\textbf{Step 2: Reduction to a scalar inequality.}
Let%
\begin{equation*}
W_{k}(r):=u_{k}(r)+v_{k}(r),\quad W(r):=u(r)+v(r),
\end{equation*}%
where $(u,v)$ is the limit solution corresponding to $(\alpha ,\beta )$.
Since $u_{k}\leq W_{k}$ and $v_{k}\leq W_{k}$, and $f,g$ are nondecreasing,
we have 
\begin{align}
\Delta W_{k}& =p(r)f(u_{k},v_{k})+q(r)g(u_{k},v_{k})  \notag \\
& \leq (p(r)+q(r))\,\big(f(W_{k},W_{k})+g(W_{k},W_{k})\big),
\label{eq:scalar-ineq}
\end{align}%
by (H2).

\medskip\noindent\textbf{Step 3: Applying the Keller--Osserman transform.}
Applying $H$ to $W_{k}$ and using $H^{\prime \prime }\geq 0$, the chain rule
gives%
\begin{equation}
\Delta H(W_{k})=H^{\prime }(W_{k})\,\Delta W_{k}+H^{\prime \prime
}(W_{k})|\nabla W_{k}|^{2}\geq H^{\prime }(W_{k})\,\Delta W_{k}.
\label{eq:chain-rule}
\end{equation}%
\emph{Justification of the chain rule:} Since $W_k \in C^2$ and $H \in C^2((0,\infty))$, the composition $H \circ W_k$ is $C^2$ wherever $W_k > 0$.
The Laplacian of the composition is
\begin{equation*}
\Delta H(W_k) = H'(W_k) \Delta W_k + H''(W_k) |\nabla W_k|^2,
\end{equation*}
which is the standard chain rule for the Laplacian. Since $H''(W_k) \geq 0$ by our earlier computation, we obtain \eqref{eq:chain-rule}.

Since $H^{\prime }(s)=-1/(f(s,s)+g(s,s))$, inequality \eqref{eq:scalar-ineq}
yields%
\begin{equation*}
\Delta H(W_{k})\geq -(p(r)+q(r)).
\end{equation*}

\medskip\noindent\textbf{Step 4: First radial integration.}
In radial coordinates, the above inequality becomes%
\begin{equation*}
\big(r^{n-1}(H(W_{k}))^{\prime }(r)\big)^{\prime }\geq -r^{n-1}(p(r)+q(r)).
\end{equation*}%
Integrating from $0$ to $r<R_{k}$ and using regularity at $r=0$ (which
implies $\lim_{t\downarrow 0}t^{n-1}(H(W_{k}))^{\prime }(t)=0$), we obtain%
\begin{equation*}
(H(W_{k}))^{\prime }(r)\geq -r^{1-n}\int_{0}^{r}s^{n-1}(p(s)+q(s))\,ds.
\end{equation*}

\medskip\noindent\textbf{Step 5: Second radial integration and lower bound.}
Integrating from $r$ to $R_{k}$ gives%
\begin{equation*}
H(W_{k}(R_{k}))-H(W_{k}(r))\geq -\int_{r}^{R_{k}}t^{1-n}\left(
\int_{0}^{t}s^{n-1}(p(s)+q(s))\,ds\right) dt.
\end{equation*}%
Letting $r\rightarrow R_{k}^{-}$, we have 
\begin{equation*}
W_{k}(R_{k})\rightarrow \infty \text{ and hence }H(W_{k}(R_{k}))\rightarrow
0.
\end{equation*}%
Thus%
\begin{equation*}
H(W_{k}(r))\leq \int_{r}^{R_{k}}t^{1-n}\left(
\int_{0}^{t}s^{n-1}(p(s)+q(s))\,ds\right) dt.
\end{equation*}%
Since $H$ is strictly decreasing, we invert to obtain%
\begin{equation*}
W_{k}(r)\geq H^{-1}\!\left(
\int_{r}^{R_{k}}t^{1-n}\int_{0}^{t}s^{n-1}(p(s)+q(s))\,ds\,dt\right) .
\end{equation*}%
Passing to the limit $k\rightarrow \infty $ yields%
\begin{equation*}
W(r)\geq H^{-1}\!\left(
\int_{r}^{R}t^{1-n}\int_{0}^{t}s^{n-1}(p(s)+q(s))\,ds\,dt\right) .
\end{equation*}%
If $R<\infty $, the right-hand side tends to $\infty $ as $r\rightarrow
R^{-} $, contradicting the fact that $(u,v)$ is entire. Hence $R=\infty $,
and we have established the quantitative lower bound \eqref{eq:lower-bound}.

\medskip\noindent\textbf{Step 6: Divergence of at least one component.}
From the previous steps (in particular, Step~1 of the proof of Theorem~\ref%
{lem:Gphi} and \eqref{le}), we already know that%
\begin{equation*}
W(r):=u(r)+v(r)\;\longrightarrow \;+\infty \quad \text{as }r\rightarrow
\infty .
\end{equation*}%
Suppose, for contradiction, that both $u$ and $v$ are bounded on $[0,\infty
) $. That is, there exist constants $M_{1},M_{2}>0$ such that%
\begin{equation*}
u(r)\leq M_{1},\quad v(r)\leq M_{2}\quad \forall r\geq 0.
\end{equation*}%
Then%
\begin{equation*}
W(r)\leq M_{1}+M_{2}\quad \forall r\geq 0,
\end{equation*}%
which contradicts $W(r)\rightarrow +\infty $. Therefore, at least one of the
two components must diverge:%
\begin{equation*}
\lim_{r\rightarrow \infty }u(r)=+\infty \quad \text{or}\quad
\lim_{r\rightarrow \infty }v(r)=+\infty .
\end{equation*}

This next figure shows 
\begin{equation*}
u\left( r\right) \text{ and }v\left( r\right) \text{ vs. }r,
\end{equation*}
with the \textquotedblleft both bounded\textquotedblright\ scenario
impossible, and the Keller--Osserman case forcing both curves upward. 
\begin{figure}[h]
\centering
\begin{tikzpicture}[scale=1.0,>=Stealth]
  % Axes
  \draw[->] (0,0) -- (8,0) node[below right] {$r$};
  \draw[->] (0,0) -- (0,5) node[left] {$u(r),\,v(r)$};

  % Labels for infinity
  \node at (-0.4,4.8) {$+\infty$};

  % Case 1: Without KO condition (only one diverges)
  \draw[thick,blue,domain=0:7,smooth,variable=\x] 
    plot ({\x},{1.5+0.5*tanh((\x-1))}) node[right] {$u(r)$ bounded};
  \draw[thick,red,domain=0:7,smooth,variable=\x] 
    plot ({\x},{0.5+0.6*\x}) node[right] {$v(r)\to+\infty$};

  % Case 2: With KO condition (both diverge) - dashed curves
  \draw[thick,blue,dashed,domain=0:7,smooth,variable=\x] 
    plot ({\x},{1.0+0.5*\x}) node[above right] {$u(r)\to+\infty$};
  \draw[thick,red,dashed,domain=0:7,smooth,variable=\x] 
    plot ({\x},{0.8+0.55*\x}) node[below right] {$v(r)\to+\infty$};

  % Legend
  \begin{scope}[shift={(4.5,1.0)}]
    \draw[blue,thick] (0,0) -- (0.6,0) node[right,black] {\small $u$ without KO};
    \draw[red,thick] (0,-0.4) -- (0.6,-0.4) node[right,black] {\small $v$ without KO};
    \draw[blue,thick,dashed] (0,-0.8) -- (0.6,-0.8) node[right,black] {\small $u$ with KO};
    \draw[red,thick,dashed] (0,-1.2) -- (0.6,-1.2) node[right,black] {\small $v$ with KO};
  \end{scope}

  % Annotation for W(r)
  \node[align=center] at (2.5,4.3) {$W(r)=u(r)+v(r)$\\$\longrightarrow +\infty$};

\end{tikzpicture}
\caption{ $W(r)\rightarrow +\infty $ forces at least one component to
diverge; under the Keller--Osserman condition both diverge. Solid lines:
without KO; dashed lines: with KO.}
\end{figure}
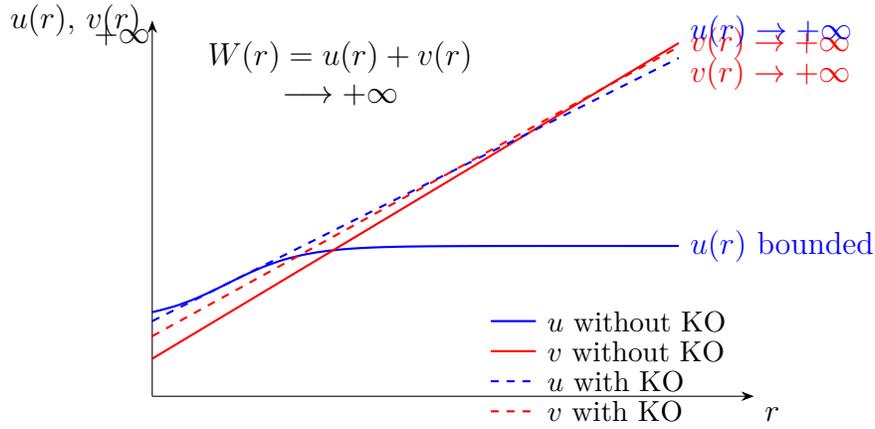

\medskip\noindent\textbf{Both components diverge under the supplementary
coupled Keller--Osserman hypothesis.}
The diagonal integrals \eqref{kelleroser} specify the superlinear growth
class, while \eqref{eq:balance} supplies the coupling needed to transfer
growth from one equation to the other.  We give two complete proofs.

\smallskip\noindent
\emph{First proof (radial flux comparison).}
Since $u,v$ are nondecreasing and their central values are strictly larger
than $\delta\ge\eta$,
(H4) applies along the whole trajectory.  Define
\[
\mathcal F_u(r):=r^{n-1}u'(r),\qquad
\mathcal F_v(r):=r^{n-1}v'(r).
\]
Regularity at the origin and the radial equations give
\[
\mathcal F_u(r)=\int_0^r s^{n-1}p(s)f(u(s),v(s))\,ds,\qquad
\mathcal F_v(r)=\int_0^r s^{n-1}q(s)g(u(s),v(s))\,ds .
\]
The integrands are nonnegative, and \eqref{eq:balance} therefore implies
\[
\kappa\mathcal F_v(r)\le\mathcal F_u(r)
\le\kappa^{-1}\mathcal F_v(r),\qquad r\ge0.
\]
After division by $r^{n-1}$ and integration from $0$ to $r$, we obtain the
pointwise two-sided estimate
\begin{equation}
\kappa\,[v(r)-\beta]\le u(r)-\alpha
\le\kappa^{-1}[v(r)-\beta].
\label{eq:component-comparison}
\end{equation}
Step~6 established $u(r)+v(r)\to\infty$.  If one component were bounded,
\eqref{eq:component-comparison} would bound the other, a contradiction.
Thus $u(r)\to\infty$ and $v(r)\to\infty$.

\smallskip\noindent
\emph{Second proof (positive Newtonian kernels).}
Reversing the order of integration in \eqref{eq:int-u} gives
\[
u(r)-\alpha=\int_0^r K_r(s)p(s)f(u(s),v(s))\,ds,
\]
and analogously for $v$, where
\[
K_r(s)=s^{n-1}\int_s^r t^{1-n}\,dt
=\frac{s}{n-2}\left[1-\left(\frac{s}{r}\right)^{n-2}\right]\ge0
\quad(0\le s\le r).
\]
Multiplication of \eqref{eq:balance} by the positive kernel $K_r$ and
integration immediately reproduces \eqref{eq:component-comparison}.  Since
$u+v$ is large, the same contradiction proves componentwise blow-up.  This
argument is independent of differentiating or comparing radial fluxes and
also identifies positivity of the Newtonian kernel as the transfer
mechanism.

\medskip\noindent\textbf{Step 7: Conclusion.}
We have now completed the construction of a positive, radial, entire solution%
\begin{equation*}
(u,v)\in C^{2}([0,\infty ))\times C^{2}([0,\infty ))
\end{equation*}%
to the system \eqref{eq:system} under the standing assumptions (H1)--(H3).
From the previous steps, we have established that the sum%
\begin{equation*}
W(r):=u(r)+v(r)
\end{equation*}%
satisfies%
\begin{equation*}
W(r)\longrightarrow +\infty \quad \text{as }r\rightarrow \infty ,
\end{equation*}%
and moreover obeys the quantitative lower bound \eqref{eq:lower-bound}
derived via the Keller--Osserman transform.

We have also shown that:

\begin{itemize}
\item Without any additional conditions added to (H1)-(H3), at least one of
the two components $u$ or $v$ must diverge to $+\infty $ as $r\rightarrow
\infty $.

\item Under the supplementary coupled Keller--Osserman hypothesis (H4),
both components necessarily diverge:%
\begin{equation*}
\lim_{r\rightarrow \infty }u(r)=\lim_{r\rightarrow \infty }v(r)=+\infty .
\end{equation*}
\end{itemize}

Thus, the constructed pair $(u,v)$ is a positive, radial, entire solution of %
\eqref{eq:system} with the stated asymptotic properties, completing the
proof of Theorem~\ref{thm:main}.
\end{proof}
\begin{remark}
In the supercritical regime, we assume that the nonlinearities admit the normalized
power-type bounds
\[
M := \max\left\{
\sup_{u>0,\;v>0}\frac{f(u,v)}{u^{a}v^{b}},
\;
\sup_{u>0,\;v>0}\frac{g(u,v)}{u^{c}v^{d}}
\right\} < \infty,
\qquad
m := \min\left\{
\inf_{u>0,\;v>0}\frac{f(u,v)}{u^{a}v^{b}},
\;
\inf_{u>0,\;v>0}\frac{g(u,v)}{u^{c}v^{d}}
\right\} > 0,
\]
together with the supercritical condition
\[
(a+b)(c+d) > 1,
\]
where a, b, c and d are positive parameters. In addition, we impose the structural condition
\[
a=c,\qquad b=d,
\]
which ensures that both reaction channels possess the same power-type profile.
Without this structural assumption, the ratio
\[
\frac{p(r)f(u,v)}{q(r)g(u,v)}
\]
cannot be uniformly controlled, and therefore no comparison principle is available
to guarantee that both components diverge to $+\infty$.

The weights are further required to satisfy the uniform comparability
\[
c_{1}\,q(r) \;\le\; p(r) \;\le\; c_{2}\,q(r),
\qquad r\ge 0,
\]
for some constants $c_{1},c_{2}>0$, while retaining the structural assumptions of (H1).

Under these hypotheses, the Keller--Osserman envelope (H3) ensures that the scalar
quantity $W=u+v$ is large, and the above comparability of the weighted reaction
channels yields the corresponding componentwise blow-up,
\[
\lim_{r\to\infty}u(r)=\lim_{r\to\infty}v(r)=\infty.
\]
\end{remark}

\textbf{Proof.}
Under the stated hypotheses we first recall that the Keller--Osserman envelope (H3)
implies that the scalar quantity
\[
W(r)=u(r)+v(r)
\]
is large, i.e.\ $W(r)\to\infty$ as $r\to\infty$. Since $u$ and $v$ are radial and satisfy
the system (1.1), their radial fluxes are given by
\[
F_u(r)=r^{\,n-1}u'(r)=\int_0^r s^{\,n-1}p(s)f(u(s),v(s))\,ds,
\qquad
F_v(r)=r^{\,n-1}v'(r)=\int_0^r s^{\,n-1}q(s)g(u(s),v(s))\,ds.
\tag{1}
\]
The normalization assumptions on $f$ and $g$ yield
\[
m\,u^{a}v^{b}\le f(u,v)\le M\,u^{a}v^{b},
\qquad
m\,u^{c}v^{d}\le g(u,v)\le M\,u^{c}v^{d},
\tag{2}
\]
for all $u,v>0$. Together with the comparability of the weights,
\[
c_{1}q(r)\le p(r)\le c_{2}q(r),\qquad r\ge 0,
\tag{3}
\]
we obtain the pointwise bounds
\[
c_{1}m\,q(r)u^{a}v^{b}
\;\le\;
p(r)f(u,v)
\;\le\;
c_{2}M\,q(r)u^{a}v^{b},
\tag{4}
\]
and similarly
\[
m\,q(r)u^{c}v^{d}
\;\le\;
q(r)g(u,v)
\;\le\;
M\,q(r)u^{c}v^{d}.
\tag{5}
\]
At this point we impose the additional structural condition
\[
a=c,\qquad b=d,
\tag{6}
\]
which expresses that both reaction channels have the same power-type profile.
Then (2), (4), and (5) simplify to
\[
m\,u^{a}v^{b}\le f(u,v)\le M\,u^{a}v^{b},
\qquad
m\,u^{a}v^{b}\le g(u,v)\le M\,u^{a}v^{b},
\]
and
\[
c_{1}m\,q(r)u^{a}v^{b}
\;\le\;
p(r)f(u,v)
\;\le\;
c_{2}M\,q(r)u^{a}v^{b},
\]

\[
m\,q(r)u^{a}v^{b}
\;\le\;
q(r)g(u,v)
\;\le\;
M\,q(r)u^{a}v^{b}.
\]
Dividing the upper and lower bounds for $p(r)f(u,v)$ by those for $q(r)g(u,v)$
we obtain constants
\[
K_{1}=\frac{c_{1}m}{M},\qquad K_{2}=\frac{c_{2}M}{m},
\]
such that
\[
K_{1}\,q(r)g(u,v)\;\le\;p(r)f(u,v)\;\le\;K_{2}\,q(r)g(u,v),
\qquad r\ge 0.
\tag{7}
\]
Integrating (7) in the flux representation (1) yields
\[
K_{1}F_v(r)\;\le\;F_u(r)\;\le\;K_{2}F_v(r),
\qquad r\ge 0.
\tag{8}
\]
Dividing by $r^{\,n-1}$ and integrating from $0$ to $r$ gives the pointwise comparison
\[
K_{1}\,[v(r)-v(0)]
\;\le\;
u(r)-u(0)
\;\le\;
K_{2}\,[v(r)-v(0)].
\tag{9}
\]
Since $W(r)=u(r)+v(r)\to\infty$ by (H3), neither $u$ nor $v$ can remain bounded;
otherwise (9) would force both components to be bounded, contradicting the
Keller--Osserman growth. Therefore,
\[
\lim_{r\to\infty}u(r)=\lim_{r\to\infty}v(r)=\infty.
\]
This completes the proof.
\qed
\section{Physical motivation: derivation of a catalytic reactor model}\label{aplicatie}

\subsection{Continuum setting and constitutive assumptions}
Consider an isothermal porous catalyst containing two dissolved or gaseous
reactants A and B.  At the pellet scale the pore geometry is homogenized:
$u(x,t)$ and $v(x,t)$ are intrinsic phase-averaged concentrations, and
\[
D_{A,\mathrm{eff}}=\frac{\varepsilon}{\tau}D_{A,m},\qquad
D_{B,\mathrm{eff}}=\frac{\varepsilon}{\tau}D_{B,m}
\]
are representative effective diffusivities, with porosity $\varepsilon$,
tortuosity $\tau$, and molecular diffusivities $D_{i,m}$; more elaborate
closures may replace these constants \cite{Aris1975,Fogler2016}.  We neglect
convection, pore-scale transients, thermal gradients, and volume change.
These assumptions identify an isothermal diffusion--reaction regime rather
than a complete reactor model.

Let $\dot\omega_A,\dot\omega_B$ denote net volumetric production rates.
Species conservation and Fick's law are
\[
\partial_t u+\nabla\!\cdot J_A=\dot\omega_A,\quad
J_A=-D_{A,\mathrm{eff}}\nabla u,
\qquad
\partial_t v+\nabla\!\cdot J_B=\dot\omega_B,\quad
J_B=-D_{B,\mathrm{eff}}\nabla v.
\]
For irreversible consumption write
$\dot\omega_A=-a_A(x)\mathcal R_A(u,v)$ and
$\dot\omega_B=-a_B(x)\mathcal R_B(u,v)$, where $a_i$ are local catalyst
activities.  At steady state one obtains, with no sign ambiguity,
\begin{equation}
\Delta u=\frac{a_A(x)}{D_{A,\mathrm{eff}}}\mathcal R_A(u,v),\qquad
\Delta v=\frac{a_B(x)}{D_{B,\mathrm{eff}}}\mathcal R_B(u,v).
\label{eq:dimensional-reactor}
\end{equation}
Thus the positive right-hand sides in \eqref{eq:system} represent
consumption rates: concentration is subharmonic inside a pellet supplied
from its exterior.  Taking
\[
p(x)=a_A(x)/D_{A,\mathrm{eff}},\quad
q(x)=a_B(x)/D_{B,\mathrm{eff}},\quad
f=\mathcal R_A,\quad g=\mathcal R_B
\]
gives exactly the mathematical system.
\subsection{Kinetics, cooperation, and dimensional consistency}

A useful phenomenological example is
\[
\mathcal R_A(u,v)=k_Au^{a+1}v^b,\qquad
\mathcal R_B(u,v)=k_Bu^cv^{d+1},
\qquad a,b,c,d\ge0.
\]
The constants $k_i$ carry the units required by the reaction orders.
Positive cross-orders describe promotion of each consumption channel by the
other reactant.  Consequently the PDE nonlinearities are cooperative in the
order-theoretic sense, although both chemical species are consumed.  More
realistic Langmuir--Hinshelwood laws can also be used when they are
monotone in the concentration range of interest; substrate inhibition,
competitive adsorption, and Arrhenius thermal feedback generally violate
(H2) and require a different analysis.

For $s,t\ge\eta$, the power-law sum admits a lower envelope
$C(s+t)^{m_*}$ for large $s+t$, where
\[
m_*=\min\{\max(a+1,c),\max(b,d+1)\}.
\]
Hence $m_*>1$ is a transparent sufficient condition for (H3).  The balanced
channel condition (H4) is stronger: it models two reaction pathways whose
weighted local rates remain comparable.  It holds, for example, for a
common kinetic factor
$\mathcal R_A=c_A\mathcal R$, $\mathcal R_B=c_B\mathcal R$ when the activity
ratio $a_A/a_B$ is bounded above and below.  This physical interpretation
also explains why diagonal superlinearity alone cannot transfer blow-up
between species.

\subsection{Spherical reduction and nondimensional groups}
For a spherical pellet of radius $L$, put $\rho=r/L$,
$U=u/C_A^*$, and $V=v/C_B^*$.  If the activities are radial,
\eqref{eq:system} becomes
\begin{align}
U_{\rho\rho}+\frac{n-1}{\rho}U_\rho
 &=\Phi_A^2\,\widehat a_A(\rho)\widehat f(U,V),\\
V_{\rho\rho}+\frac{n-1}{\rho}V_\rho
 &=\Phi_B^2\,\widehat a_B(\rho)\widehat g(U,V),
\end{align}
where the squared generalized Thiele moduli have the form
\[
\Phi_A^2=\frac{L^2R_A^*}{D_{A,\mathrm{eff}}C_A^*},\qquad
\Phi_B^2=\frac{L^2R_B^*}{D_{B,\mathrm{eff}}C_B^*}.
\]
They compare intrinsic reaction and intrapellet diffusion time scales.
Regularity imposes $U_\rho(0)=V_\rho(0)=0$.  At the physical pellet surface
one ordinarily prescribes either concentrations or external mass-transfer
conditions
\[
D_{A,\mathrm{eff}}u_r(L)=k_{m,A}(u_b-u(L)),\qquad
D_{B,\mathrm{eff}}v_r(L)=k_{m,B}(v_b-v(L)).
\]
The associated Biot numbers $k_{m,i}L/D_{i,\mathrm{eff}}$ distinguish
external from intrapellet transport resistance; effectiveness factors then
measure the volume-averaged rate relative to the surface-rate benchmark
\cite{Aris1975,Fogler2016}.

\subsection{Meaning of the whole-space idealization}
The domain $\mathbb R^n$ is not a literal finite pellet.  It is an
idealization of a radially heterogeneous reactive medium, or the limit of
pellets whose outer boundary is moved outward.  Profiles such as
\[
p(r)=A(1+r)^{-\mu},\qquad q(r)=B(1+r)^{-\nu},
\qquad \mu,\nu>2,
\]
represent activity attenuation and satisfy (H1).  Entire large solutions
describe a mathematical loss of uniform concentration control under
far-field forcing.  They should not be interpreted as finite-mass steady
states of an isolated reactor: unbounded concentrations signal the failure
of the isothermal, constant-property, or infinite-reservoir approximation.
In an exothermic pellet a temperature equation and Arrhenius kinetics would
be required to model thermal ignition or runaway.

\section{Conclusion \label{5}}

We established a threshold construction for positive radial entire
solutions of a weighted cooperative elliptic system.  The decisive scalar
estimate applies to $W=u+v$ and shows that a boundary central datum produces
$W(r)\to\infty$, together with the explicit tail-potential bound
\eqref{eq:lower-bound}.  Total blow-up does not, by itself, identify the
asymptotics of each component.  The balanced-channel assumption (H4)
resolves this issue through the quantitative comparison
\eqref{eq:component-comparison}; the flux and Newtonian-kernel proofs expose
two complementary forms of the same positivity mechanism.

The distinction between scalar Keller--Osserman growth and genuine coupling
is one of the main structural conclusions.  It prevents an invalid
componentwise inference in decoupled systems and gives assumptions that can
be checked directly in applications.  The reactor derivation further shows
how finite Newtonian potentials correspond to decaying catalyst activity
and why whole-space large profiles must be interpreted as an asymptotic
idealization.

Natural extensions include sharp necessary conditions for balance,
component-specific asymptotic equivalents, nonradial threshold solutions,
and nonisothermal systems with a temperature equation.  Each would require
tools beyond the radial order and positive-kernel arguments used here.

\section*{Declaration of Generative AI and AI-assisted technologies in the writing process}
The authors utilized the free AI tool Grammarly to improve the manuscript’s grammar. The Section~\ref{aplicatie} was prepared with the support of Microsoft Copilot AI in Microsoft Edge, drawing on the ideas of \cite{KO56} in the scalar case to ensure a clear and coherent integration of the applied perspective into the theoretical framework. Following its use, the authors reviewed and edited the content as necessary, taking full responsibility for the final publication.
\section*{Declaration of competing interest}
This work does not have any conflicts of interest.
\section*{Funding}
This research did not receive any specific grant from funding agencies in the public, commercial, or not-for-profit sectors.
\section*{Data availability}
No data was used for the research described in the article.
%% The Appendices part is started with the command \appendix;
%% appendix sections are then done as normal sections
\section*{Acknowledgments}
The author expresses his gratitude to the referees.

\end{document}